\documentclass[10pt,a4paper,twoside,leqno]{article}
\usepackage{amsfonts,amssymb}
\usepackage{eufrak}
\usepackage{amscd}

\font\Bigtit=cmr10 scaled \magstep 4
\font\ebf=cmbx8
\font\erm=cmr8

\usepackage{graphicx}

\begin{document}

\thispagestyle{empty}

\begin{flushright}
PL ISSN 0459-6854
\end{flushright}
\vspace{0.5cm}
\centerline{\Bigtit B U L L E T I N}
\vspace{0.5cm}
\centerline{DE \ \  LA \ \  SOCI\'ET\'E \ \  DES \ \  SCIENCES \ \ ET \ \ DES \
\ \ LETTRES \ \ DE \ \ \L \'OD\'Z}
\vspace{0.3cm}
\noindent 2003\hfill Vol. LIII
\vspace{0.3cm}
\hrule
\vspace{5pt}
\noindent Recherches sur les d\'eformations \hfill Vol. XLII
\vspace{5pt}
\hrule
\vspace{0.3cm}
\noindent pp.~11--25

\vspace{0.7cm}

\noindent {\it Andrzej K. Kwa\'sniewski}

\vspace{0.5cm}

\noindent {\bf ON DUALITY TRIADS}

\vspace{0.5cm}

\noindent {\ebf Summary}

{\small
In this paper we introduce a   duality triads`  notion. These are dual
recurrences \linebreak[4] [1, 2] as used in dynamical data bases theory  completed by a
third  pertinacious  relation.
After duality triads  being defined  several representative examples of
them are given.  One illustrates these examples with help of Pascal-like
triangles which as a matter of fact contain in the illustrative way the
clue information on the triad system.
$q$-Gaussian triads as well as  Fibonomial   triangle [3--8]  and duality
triads in finite operator calculus [9--17]  are presented in the paper that
follows this one.}

\vspace{0.5cm}

\renewcommand{\thesubsection}{\arabic{subsection}.}

\subsection{The notion of duality triads' system}

Let us consider a kind of random walking process on   $N \cup \{0\}$ set of
nonnegative integers. Let  $c_{n,k}$ be the number {\it of ways to reach level
$k$ in $n$ steps starting from the  level}  0. Then  [1, 2] :
\begin{eqnarray}
&&c _{n+1,k}=i _{k-1}c _{n,k-1}+q_{k}c_{n,k}+d_{k+1}c
_{n,k+1}, \label{eq*0} \\
&&\nonumber \\
&&c _{0,0}=1;\ c _{0,k}=0\ {\rm
for}\ k>0. \nonumber  \end{eqnarray}

\noindent Here numbers $i_{k},\ q_{k},\ d_{k};\ k \geq 0$ -- independent of
``discrete time'' parameter  $n$   are defined as:
$$
  i_{k}=N{\rm pos}(I,k),\quad q_{k}=N{\rm pos}(Q,k),\quad d_{k}=N{\rm pos}(D,k);\quad k \geq
0, $$

\noindent where   $N{\rm pos}(O,k) =$ {\it number of possibilities to perform the
operation O on the file of size $k$.}
One of the most important enumerative object of the analysis of data
structures [1, 2] is the sequence $\{ H_{k,l,n}\}$, where  $H_{k,l,n} =$
{\it number of histories of length n with initial height k and final height
l}. While starting with an empty file it is natural then to consider in
this special situation  the above introduced array  of nonnegative integer
numbers  i.e. -- say it again --

\vspace{3mm}

\noindent $c _{n,k}=H_{0,k,n} =$ {\it number of ways to reach level k in n
steps starting from the level} 0.

\vspace{3mm}

\noindent With this in mind one recognizes  that for the one-step
transition from time $n$ to $n+1$ the array of numbers $c_{n,k}$   must be
the unique  solution of  the   recurrence Eq.~(\ref{eq*0}).

Let now $\{ \Phi _{n} \}_{n \geq 0}$  be a polynomial sequence  (i.e.
${\rm deg} \Phi_{n}(x) = n$) determined by the following   recurrence equation
``dual with
respect to do (\ref{eq*0})'':
\begin{eqnarray}
&& x \Phi _{n}(x)=d _{n}\Phi_{n-1}(x)+q_{n}\Phi_{n}(x)+i_{n}\Phi_{n+1}(x),
\label{eq**0}\\
&&\nonumber \\
&&\Phi _{0}(x)=1,\quad \Phi_{-1}(x)=0;\quad n \geq 0.\nonumber
\end{eqnarray}

\noindent In combinatorics one calls   the  recurrence Eq. (\ref{eq**0})
-- {\it  the dual recurrence} with respect to the recurrence Eq.~(\ref{eq*0})
because from Eq.~(\ref{eq*0})    equivalent to Eq.~(\ref{eq**0})
one derives the so-called duality relations between
polynomial solutions of Eq.~(\ref{eq**0}) and monomials $x^{n}$. Namely, one
may see that (with $i_{k}, q_{k}, d_{k}; k \geq 0$ -- independent of
``discrete time'' parameter  $n$):

\vspace{3mm}

\noindent {\bf Lemma 1.}
\begin{equation}
x_{n}=\sum\limits_{k\geq 0}c_{n,k}\Phi _{k}(x) \quad n \geq 0.\label{eq***0}
\end{equation}

\vspace{3mm}

\noindent {\it Proof.} The identity (\ref{eq***0}) is trivial for  $n = 0$.
Let it be then true for $n\geq 0$. Then
{\setlength\arraycolsep{1pt}
\begin{eqnarray*}
x^{n+1}&=& x \sum\limits_{k \geq 0} C _{n,k}\phi _{k}(x)=\{ {\rm see\
(3)}\}=\sum\limits_{k\geq 0} C _{n,k} \left\{ d_{k}\Phi_{k-1}(x)+q_{k}\Phi
_{k}(x)+i_{k}\Phi _{k+1}(x) \right\}=;\\
&=&\left\{ {\rm due\ to}\ \sum\limits _{k\geq
0}C_{n,k}i_{k}\Phi_{k+1}=\sum\limits _{k=0}C_{n,k-1}i_{k-1} \Phi _{k}\{i
_{-1}=0 \}\ {\rm and\ due\ to}\right.\\
&&\left. \sum\limits_{k\geq 0} C _{n,k}d_{k}\Phi
_{k-1} = \sum\limits_{k \geq 1}C_{n,k}d_{k}\Phi _{k-1}\{ d_{0}=0
\}=\sum\limits_{k\geq 0}C_{n,k+1}d_{k+1}\Phi _{k} \right\} =\\
&=& \sum\limits _{k\geq 0}\left\{d_{k+1}
C_{n,k+1}+q_{k}C_{n,k}+i_{k-1}C_{n,k-1} \right\} \Phi _{k} = \sum\limits
_{k \geq 0} C _{n+1,k} \phi _{k}(x). \end{eqnarray*}}

\vspace{3mm}

\noindent {\it Definition 1}. $\{ \Phi_{n}  \}_{n \geq 0 }$  are then said
to be {\it triad polynomials.}

\vspace{3mm}

The formula (\ref{eq***0}) allows us  to interpret the counting
sequence $c_{n,k}$ as the sequence of expansion coefficients of monomials
$x_{n}$ in the basis of  polynomials $\{ \Phi _{n} \}_{n\geq 0}$ as in the
framework of umbral or finite operator calculus. The   array of  numbers
$c_{n,k}$  is interpreted as special connection constants` lower
triangle infinite array [10--12] and other triangles. (The Fibonomial
[3--8] triangle no triad case is treated separately at the end).

Let us note that restoring either algorithms or formulas for connection
constants  is one of the central problems in finite operator calculus in its
classical formulation  [11, 12] or in its later on  extended formulation
(see  [13--17]  and references therein).

\subsubsection*{Remarks and Information I }

Note however  that  from (\ref{eq***0})  it does not follow  that numbers $c_{n,k}$
 satisfy a  second order recurrence  (\ref{eq*0}). Compare
$$
x^{n}=\sum_{k=0}^{n}\left(\begin{array}{c} n\\k \end{array}
\right)k^{n-k}A_{k}(x),\quad n \geq 0$$

\noindent where $\{ A_{n} \}_{n \geq 0}$ represents the Abel
binomial polynomial sequence $A_{n}(x)=[x(x+n)^{n-1}],\ n \geq 0$ (the
recurrence for $\{ A_{n} \}_{n \geq 0}$ see: p. 73 in [ 16]). Neither
it follows that corresponding
polynomial sequence (i.e. ${\rm deg} \Phi _{k}(x) = k$)  exists at all.
Compare with Euler numbers $\langle\begin{array}{c} n\\k \end{array}\rangle$ [18] for which
$$
x^{n}=\sum\limits_{k}\left\langle \begin{array}{c} n\\k \end{array}
\right\rangle \left( \begin{array}{c} x+k\\ n \end{array} \right) \neq k
\quad {\rm for}\quad k\neq n,$$

\noindent for numbers $i_{k},
q_{k}, d_{k}; k \neq 0$ -- dependent of ``discrete time'' parameter $n$.
$\{ A_{n} \}_{n \geq 0}$ see: p. 73 in [16]). Neither it follows that
corresponding polynomial sequence (i.e. ${\rm deg} \Phi_{k}(x) = k$) exists
at all for numbers $i_{k}, q_{k}, d_{k}; k \geq 0$ -- dependent of
``discrete time'' parameter $n$.

Compare with Euler numbers $\langle \begin{array}{c} n\\k \end{array}
\rangle$ [18] for which
$$
x^{n}=\sum\limits_{k}\left\langle \begin{array}{c} n\\k \end{array}
\right\rangle \left( \begin{array}{c} x+k\\ n \end{array} \right) \neq k
\quad {\rm for}\quad k\neq n,$$

\noindent while the corresponding
recurrence for Euler numbers $\langle \begin{array}{c} n\\k \end{array}
\rangle$ reads:
$$
\left\langle \begin{array}{c} n+1\\k \end{array} \right\rangle =(k+1)
\left\langle \begin{array}{c} n\\k \end{array} \right\rangle +
(n+1-k) \left\langle \begin{array}{c} n\\k-1 \end{array} \right\rangle ,$$
$$k>0;\
\left\langle \begin{array}{c} 0\\0 \end{array} \right\rangle =1,\
\left\langle \begin{array}{c} 0\\k \end{array} \right\rangle =0,\
k>0.
 $$

\noindent There are many  polynomial sequences of distinguished importance
which are not solutions of recurrence of the type (\ref{eq**0}),  for
example Euler or    (see: p. 203 in [11])  exponential polynomials are not
triad polynomial sequences.

Existence of dual triad is then a quite specifically qualified property of
corresponding objects.    Let then  ${\rm deg} \Phi _{k}(x) = k.$

\vspace{3mm}

\noindent {\it Definition 2}.   The ensemble of   (\ref{eq*0}),
(\ref{eq**0}), (\ref{eq***0}) shall be called  the duality  triad system
or just {\it duality triad} -- for short. Polynomials $\{ \Phi _{n} \}_{n \geq
0}$ are then said to be {\it triad polynomials}.

\vspace{3mm}

To this end let us notice that the recurrences  (\ref{eq*0}) and
(\ref{eq**0}) might be interpreted as  defining  simple discrete-time   semi-group
dynamical  system under iteration of  the ``one step transition matrix'' $X$:
\begin{equation}\label{eq1}
X=\left(\begin{array}{cccccc} q_{0}&i_{0}&0&0&0&\ldots\\
d_{1}&q_{1}&i_{1}&0&0&\ldots\\
0&d_{2}&q_{2}&i_{2}&0&\ldots\\
0&0&d_{3}&q_{3}&i_{3}&\ldots\\
\ldots&\ldots&\ldots&\ldots&\ldots&\ldots \end{array}
\right). \end{equation}

This  transition matrix   $X$   may be interpreted as the matrix of the
difference equation (\ref{eq**0}),  i.e.   (\ref{eq**0})  may be rewritten
in a matrix form and then be looked upon as the eigenvector equation:
\begin{equation}\label{eq2}
\left(\begin{array}{cccccc} q_{0}&i_{0}&0&0&0&\ldots\\
d_{1}&q_{1}&i_{1}&0&0&\ldots\\
0&d_{2}&q_{2}&i_{2}&0&\ldots\\
0&0&d_{3}&q_{3}&i_{3}&\ldots\\
\ldots&\ldots&\ldots&\ldots&\ldots&\ldots \end{array}
\right)
\left( \begin{array}{c} \Phi_{0}\\\Phi_{1}\\\Phi_{2}\\\Phi_{3}\\ \ldots
\end{array} \right)
=x \left( \begin{array}{c} \Phi_{0}\\\Phi_{1}\\\Phi_{2}\\\Phi_{3}\\ \ldots
\end{array} \right).\end{equation}

\noindent One also easily recognizes   that  for the one-step transition
from time (or level)  $n$ to $n+1$:
$$
c _{n+1,k}=i _{k-1}c _{n,k-1}+q_{k}c_{n,k}+d_{k+1}c
_{n,k+1} $$

\noindent we have
$$c _{0,0}=1;\quad c _{0,k}=0\quad {\rm
for}\quad k>0.   $$

\noindent so the  recurrence (\ref{eq3})  in its array form should  read as
follows:
\begin{eqnarray}
&&\left(\begin{array}{ccccc} c_{0,0}&c_{0,1}&c_{0,2}&c_{0,3}&\ldots\\
c_{1,0}&c_{1,1}&c_{1,2}&c_{1,3}&\ldots\\
c_{2,0}&c_{2,1}&c_{2,2}&c_{2,3}&\ldots\\
c_{3,0}&c_{3,1}&c_{3,2}&c_{3,3}&\ldots\\
\ldots&\ldots&\ldots&\ldots&\ldots
\end{array} \right)
\left(\begin{array}{cccccc} q_{0}&i_{0}&0&0&0&\ldots\\
d_{1}&q_{1}&i_{1}&0&0&\ldots\\
0&d_{2}&q_{2}&i_{2}&0&\ldots\\
0&0&d_{3}&q_{3}&i_{3}&\ldots\\
\ldots&\ldots&\ldots&\ldots&\ldots&\ldots \end{array}
\right) \nonumber\\
&& \label{eq3}\\
&=&
\left(\begin{array}{ccccc}
c_{1,0}&c_{1,1}&c_{1,2}&c_{1,3}&\ldots\\
c_{2,0}&c_{2,1}&c_{2,2}&c_{2,3}&\ldots\\
c_{3,0}&c_{3,1}&c_{3,2}&c_{3,3}&\ldots\\
c_{4,0}&c_{4,1}&c_{4,2}&c_{4,3}&\ldots\\
\ldots&\ldots&\ldots&\ldots&\ldots
\end{array} \right). \nonumber
\end{eqnarray}

\noindent Introducing  the row  vector $C_{n}$  (``the state of a system'')
as
$$C_{n} = (c_{n,0}, c_{n,1}, c_{n,2}, c_{n,3}, \ldots),$$

\noindent  the recurrences
(\ref{eq3}) may be written  in the form of discrete-time dynamical under
iteration system
\begin{equation}\label{eq4}
C_{n+1}=C_{n}X;\quad n\geq 0;\quad {\rm hence}\ C_{n}=C_{0}X^{n};\ n \geq 0,
\end{equation}

\noindent where $C_{0} = (1, 0, 0, 0, \ldots)$ or
\begin{equation}\label{eq4prime}
\left(\begin{array}{c} \vec{C}_{0} \\ \vec{C}_{1}\\ \vec{C}_{2}\\
\vec{C}_{3}\\ \ldots \end{array}  \right)X =
\left(\begin{array}{c} \vec{C}_{1}\\ \vec{C}_{2}\\
\vec{C}_{3}\\ \vec{C}_{4} \\ \ldots \end{array}  \right),
 \end{equation}

\noindent under the identification
$$
\left( c_{nk}  \right)=C=  \left(\begin{array}{c} \vec{C}_{0} \\ \vec{C}_{1}\\ \vec{C}_{2}\\
\vec{C}_{3}\\ \ldots \end{array}  \right), \quad \vec{C}_{n}=\left(
C_{n,0}, C_{n,1}, C_{n,2},\ldots , C_{n,k}, \ldots \right). $$

\noindent Here $(c_{nk})=C$ constitutes the lower triangle infinite array of  numbers
$c_{n,k}$     interpreted in the finite operator calculus  as
special connection  constants` [10--12].

Naturally  the recurrence  (\ref{eq3}) or  (\ref{eq4})  or
(\ref{eq4prime}) may be written also in the form  (with  transposed
$X^{T}$ {\it dual}   to  $X$) as compared  to (\ref{eq6}):
\begin{equation}\label{eq5}
C_{n+1}^{T}=X^{T}C_{n}^{T};\ n\geq 0\quad {\rm hence}\quad C_{n+1}^{T}=\left( X^{T}
\right)^{n}C_{0}^{T};\ n \geq 0. \end{equation}

\noindent Indeed,  under  the notation
$$
\Phi =\left(\begin{array}{c} \Phi_{0}\\ \Phi_{1}\\ \Phi_{2}\\ \Phi_{3} \\
\ldots \end{array} \right) $$

\noindent the recurrence  (\ref{eq**0})
takes the form  of an eigenvector and eigenvalue equation
\begin{equation}\label{eq6}
x\Phi = X\Phi .
\end{equation}

\noindent As for the iteration  of  the transition matrix $X$ the following
{\it interpretation} is pertinent [1]: the $kl$-th entry of $X_{n} =$
{\it number of ways of going from the level $k \to to \to l$   in $n$
steps}.

\subsection{First  Examples  of  triads }

Here come some classical and then new  examples of  triads.

\subsubsection{Pascal triad}

This is at the same time an example of one of the simplest  Appell
polynomials` [19, 20, 13, 19, 22]  sequence of triad polynomials ${\Phi
_{n}n(x)}_{n \geq 0};\ \Phi _{n}(x) =  (x-1)^{n}$.  Namely: the  choice
$i_{k} = 1 = q_{k}, d_{k} = 0; k \geq 0$ leads to    $c_{n,k}  = \left(
\begin{array}{c} n\\k \end{array} \right)$ and consequently to the Pascal
triangle and
\begin{eqnarray*}
&&c _{n+1,k}=i _{k-1}c _{n,k-1}+q_{k}c_{n,k}+d_{k+1}c
_{n,k+1},\\
&& \\
&&c _{0,0}=1;\ c _{0,k}=0\ {\rm
for}\ k>0,   \end{eqnarray*}
\begin{eqnarray*}
&& x \Phi _{n}(x)=d _{n}\Phi_{n-1}(x)+q_{n}\Phi_{n}(x)+i_{n}\Phi_{n+1}(x) ,
\\
&&\\
&&\Phi _{0}(x)=1,\quad \Phi_{-1}(x)=0;\quad n \geq 0,
\end{eqnarray*}
$$
x_{n}=\sum\limits_{k\geq 0}c_{n,k}\Phi _{k}(x) \quad n \geq 0.
$$

\noindent The Pascal triangle might be obtained  while considering subsequent powers of
$$
X =\left( \begin{array}{cccccc} 1&1&0&0&0&\ldots\\
0&1&1&0&0&\ldots\\
0&0&1&1&0&\ldots\\
0&0&0&1&1&\ldots\\
\ldots&\ldots&\ldots&\ldots&\ldots&\ldots  \end{array} \right)
$$

\noindent and the calculating  $C_{n} = C_{0} X^{n}$;   $n \geq 0$.  An
easy generalization is obtained for  connection constants  $c_{n,k}  =
\left( \begin{array}{c} n\\k \end{array} \right)s^{n-k}$, where $s$  is a
parameter.

\subsubsection{Stirling triad}

The choice  $i_{k} = 1, q_{k} = k; k \geq 0, d_{k} = 0; k>0$  leads to
$c_{n,k} = \left( \begin{array}{c} n\\k \end{array} \right)$  i.e. Stirling
numbers of the second type for which we have (3):
$$
x_{n}=\sum\limits_{k\geq 0}c_{n,k}\Phi _{k}(x), \quad n \geq 0.
$$

\noindent where  triad polynomials are $\Phi_{k}(x) = x^{\underline{k}}=
x(x-1)(x-2)\ldots(x-k+1)$ and consequently, following (2) and (1),
\begin{eqnarray*}
&& x \Phi _{n}(x)=d _{n}\Phi_{n-1}(x)+q_{n}\Phi_{n}(x)+i_{n}\Phi_{n+1}(x),
\\
&&\\
&&\Phi _{0}(x)=1,\quad \Phi_{-1}(x)=0;\quad n \geq 0,
\end{eqnarray*}
\begin{eqnarray*}
&&c _{n+1,k}=i _{k-1}c _{n,k-1}+q_{k}c_{n,k}+d_{k+1}c
_{n,k+1},\\
&& \\
&&c _{0,0}=1;\ c _{0,k}=0\ {\rm
for}\ k>0 .  \end{eqnarray*}

\noindent The above can be illustrated by the corresponding well known
II-type -- Stirling  triangle
$$
\begin{array}{ccccccccccccc}
&&&&&&{\bf 1}&&&&&&\\
&&&&&{\bf 0}&&{\bf 1}&&&&&\\
&&&&{\bf 0}&&1&&{\bf 1}&&&&\\
&&&{\bf 0}&&1&&3&&{\bf 1}&&&\\
&&{\bf 0}&&{\bf 1}&&7&&6&&{\bf 1}&&\\
&{\bf 0}&&1&&{\bf 15}&&25&&10&&{\bf 1}&\\
\bullet&&\bullet&&\bullet&&\bullet&&\bullet&&\bullet&&\bullet \end{array}
$$

\noindent which supplies in a compact  form the full defining information
on the triad.  The II-Stirling triangle might be obtained  while
considering subsequent powers of
$$
X =\left( \begin{array}{cccccc} 0&1&0&0&0&\ldots\\
0&1&1&0&0&\ldots\\
0&0&2&1&0&\ldots\\
0&0&0&3&1&\ldots\\
\ldots&\ldots&\ldots&\ldots&\ldots&\ldots  \end{array} \right)
$$

\noindent and then calculating  its  subsequent  rows  $C_{n} = C_{0}
X^{n};\ n \geq 0$.

Persistent roots polynomials  $\Phi _{k}(x)=x^{\underline{k}},\ k \geq 0$,
constitute a classical example  of polynomials of binomial type  which  are
basic polynomials for  the  forward difference calculus delta operator
$\Delta \equiv E -id$, where
$$E^{a}\equiv \exp \left\{a \frac{d}{dx}  \right\}; \quad a \in
\mathbb{R}$$

\noindent denotes a {\it shift operator}, i.e.
$$
\left( E^{a}f  \right)(x)=f\left( x+a \right);\quad {\rm for} \quad f:R \to
R. $$

\subsubsection*{Remarks and Information  II}

Recall: $(c_{nk}) = C$ constitutes the lower triangle infinite array of numbers
$c_{n,k}$     interpreted in finite operator calculus  as special
connection  constants [10--12],  [23--25]. In the  general  case of  two
arbitrary $\{q_{n}\}_{n \geq 0}$  and  $\{p_{n} \}_{n \geq 0}$
{\it polynomial sequences} i.e. ${\rm deg}\ p_{n} = {\rm deg}\ q_{n} = n,\ n
\geq 0 $; connection constants $(L_{n,k})$ are defined via
$p_{n}(x)=\sum_{0\leq k \leq n}L_{n,k}q_{k}(x)$, $n \geq 0$ -- using the
notation of [23, 24]. In [23] one derives explicit formulas for connection
constants $(L_{n,k})$ in terms of roots of {\it monic}
polynomial sequences $\{q_{n}\}_{n \geq 0}$  and  $\{p_{n} \}_{n \geq 0}$.
In the special case $\{q_{n}\}_{n \geq 0}$  of being any monic polynomial
sequence with persistent roots [26] one states (see: Propositions 10 and 9
in [23]) that $i_{k} = 1,\ q_{k} \neq 0;\ k\geq 0,\ d_{k} = 0;\ k>0$ (as in
Examples 1--3). Namely:
\begin{eqnarray}
&&L_{n+1,k}=L_{n,k-1}+\left( r_{k+1}-s_{n+1} \right)L_{n,k},\nonumber \\
&&\label{eqL}\\
&&L_{0,0}=1,\quad L_{0,-1}=0,\quad n,k\geq 0, \nonumber
\end{eqnarray}

\noindent where  for
$$p_{n}(x)=\sum_{k\geq 0}L_{n,k}q_{k}(x), \quad n \geq 0,
$$

\noindent $[r] \equiv \{r_{0}, r_{1}, r_{2},\ldots , r_{k}, \ldots  \}$ is the root
sequence determining $\{ q_{n} \}_{n \geq 0}$ while $[s]\equiv \{ s_{0},
s_{1}, s_{2}, \ldots , s _{k}, \ldots \}$ is the root sequence determining
$\{ p_{n} \}_{n \geq 0}$. The authors of [23] call such connection
constants $(L_{n,k})$ -- the {\it generalized Lah numbers} because with the
choice $[r] \equiv \{0,1,2,\ldots ,k,\ldots \}$ and $[s] \equiv
\{0,-1,-2,\ldots ,-k,\ldots \}$ the Lah numbers and connection constants
$L_{n,k}$ coincide as can be seen from
$$
X^{\overline{n}}=\sum\limits_{1\leq k \leq n}\left( \begin{array}{c}
n-1\\k-1\end{array} \right) \frac{n!}{k!} x ^{\underline{k}};\quad n,k\geq
0,$$

\noindent where -- recall --
$$
c_{nk}=\left(\begin{array}{c}
n-1\\n-k\end{array}  \right) \frac{n!}{k!}
$$

\noindent are Ivo Lah numbers.

  As for the treatment of the  Fibonacci and Lucas  numbers as cumulative
connection constants $K_{n}=\sum_{0\leq k \leq n}L_{n,k}$  (sums of  $n$-th
row numbers in a corresponding triangle),  see [24].

\subsubsection*{Remarks and  Information III}

The requirement  that  persistent root  polynomial sequences $\{q_{n}\}_{n
\geq 0}$ and $\{ p_{n} \}_{n \geq 0}$ should be  monic, leads  to
(\ref{eqL}). If not monic -- as it is for example in the case with Newton-Gregory
polynomials (see: (1.6) in [27]) $\Phi _{k}(x)=\left(\begin{array}{c}
x\\k\end{array} \right)$, $k \geq 0$, one gets different recurrences: here
-- for Newton-Gregory polynomials the triad is given by
\begin{eqnarray}
&&x\Phi_{k}(x)=(k+1)\Phi _{k+1}(x)+k\Phi _{k}(x),\nonumber\\
&&\label{eq**}\\
&&\Phi_{0}(x)=1,\quad \Phi_{-1}(x)=0;\quad n \geq 0, \nonumber
\end{eqnarray}

\noindent  and, consequently,
\begin{eqnarray}
&&L_{n+1,k}=kL_{n, k-1}+kL_{n,k}, \nonumber \\
&&\label{eq*}\\
&&L_{0,0}=1,\quad L_{0,-1}=0,\quad n,k \geq 0,\nonumber
\end{eqnarray}
\begin{equation}\label{eq***}
x^{n}=\sum\limits_{k\geq 0} k! \left\{ \begin{array}{c} n\\k \end{array}
\right\} \Phi_{k}(x),
\end{equation}

\noindent so the Newton-Gregory triangle   is then of the form
$$
\begin{array}{ccccccccccccc}
&&&&&&{1}&&&&&&\\
&&&&&{0}&&{1}&&&&&\\
&&&&{0}&&1&&{2}&&&&\\
&&&{0}&&1&&6&&{6}&&&\\
&&{0}&&{1}&&14&&36&&{24}&&\\
&{0}&&1&&{30}&&150&&240&&{120}&\\
\bullet&&\bullet&&\bullet&&\bullet&&\bullet&&\bullet&&\bullet \end{array}
$$

\noindent This triangle can be obtained  while considering subsequent powers of
$$
X =\left( \begin{array}{cccccc} 0&1&0&0&0&\ldots\\
0&1&2&0&0&\ldots\\
0&0&2&3&0&\ldots\\
0&0&0&3&4&\ldots\\
\ldots&\ldots&\ldots&\ldots&\ldots&\ldots  \end{array} \right)
$$

\noindent and then calculating  its  subsequent  rows $C_{n}= C_{0} X^{n}$;
$n \geq 0$.

\subsection{Stirling signed triad }

The choice $i_{k} = 1,\ q_{k} = -k$; $k \geq 0,\ d_{k} = 0, k>0$  gives
$$c_{n,k} =\left(\begin{array}{c} n\\k \end{array}  \right)(-1)^{n-k},\quad
k \geq 0,$$

\noindent where $\left(\begin{array}{c} n\\k \end{array}  \right)$ are
Stirling numbers of the second type. Thus  we have:
\begin{equation}\label{eq***2}
x^{n}=\sum\limits _{k \geq 0} \left(\begin{array}{c} n\\k \end{array}
\right) (-1)^{n-k}\Phi_{k}(x),
\end{equation}

\noindent where
$$
\Phi_{k}(x) = x^{\overline{k}} = x(x+1)(x+2)\ldots(x+k-1)$$

\noindent and, consequently,
\begin{eqnarray}
&&x \Phi_{k}(x)=\Phi_{k+1}(x)-k\Phi_{k}(x),\nonumber \\
&& \label{eq**2}\\
&&\Phi_{0}(x)=1,\quad \Phi_{-1}(x)=0;\quad n \geq 0, \nonumber
\end{eqnarray}
\begin{eqnarray}
&&c_{n+1,k}=c_{n,k-1}-kc_{n,k},\nonumber \\
&& \label{eq*2}\\
&&c_{0,0}=1,\quad c_{0,k}=0\quad {\rm for}\quad k>0. \nonumber
\end{eqnarray}

The above might be illustrated by the corresponding II-Stirling signed triangle
$$
\begin{array}{ccccccccccccc}
&&&&&&{1}&&&&&&\\
&&&&&{0}&&{1}&&&&&\\
&&&&{0}&&-1&&{1}&&&&\\
&&&{0}&&1&&-3&&{1}&&&\\
&&{0}&&{-1}&&7&&-6&&{1}&&\\
&{0}&&1&&{-15}&&25&&-10&&{1}&\\
\bullet&&\bullet&&\bullet&&\bullet&&\bullet&&\bullet&&\bullet \end{array}
$$

\noindent Persistent roots polynomials  $\Phi_{k}(x) = x^{\overline{k}},\
k \geq 0$ are strictly related  to Stirling numbers  $\left[
\begin{array}{c} n\\k \end{array} \right]$ of the first kind triangle (see the
correspondent statement and formula 6.33 in [18]). Let then note that
the persistent roots polynomial sequence  $\{
x^{\overline{k}} \}_{k \geq 0}$
constitutes simultaneously another classical example  of polynomials of
binomial type [11--17]  which are now basic polynomials for backward
difference calculus delta operator $\nabla$:
$$
\nabla \equiv id -E^{-1};\quad \left( \nabla f \right)(x)=f(x)-f(x-1)\quad
{\rm for}\quad f:R \to R. $$

\subsection{Hermite triad }

With  the following choice of    numbers:
$$
q_{k}=0;\quad k \geq 0,\quad d_{k}=k;\quad k \geq 1,\quad i_{k}=1;\quad k
\geq 0, $$

\noindent  the recurrence for monic polynomials    takes the form
\begin{eqnarray}
&&xH_{k}(x)=kH_{k-1}(x)+H_{k+1}(x),\quad k \geq 0, \nonumber\\
&& \label{eq**3} \\
&&H_{0}=1,\quad H_{-1}=0.\nonumber
\end{eqnarray}

\noindent which  is the recurrent definition of Hermite monic polynomials
(monic: i.e. the coefficient of   $x_{n}$ is one; ${\rm deg}\Phi_{n} = n$).
This is the example of nontrivial, important  Appell polynomial sequence
[19, 20, 13, 21, 22].  Naturally   Hermite classical orthogonal polynomials
 are solutions of the recurrence (\ref{eq**3}) which is  dual to the
recurrence (\ref{eq*3}) :
\begin{eqnarray}
&&c_{n+1,k}=c_{n,k-1}+(k+1)c_{n,k+1},\nonumber\\
&&\label{eq*3}\\
&&c_{0,0}=1;\quad c_{0,k}=0\quad {\rm for}\quad k>0. \nonumber
\end{eqnarray}

\noindent Hence  we also have
\begin{equation}\label{eq***3}
x^{n}=\sum\limits_{k\geq 0}c_{n,k}H_{k}(x)\quad n \geq 0.
\end{equation}

\noindent The above can be illustrated by the corresponding Hermite  triangle
$$
\begin{array}{ccccccccccccc}
&&&&&&{1}&&&&&&\\
&&&&&{0}&&{1}&&&&&\\
&&&&{1}&&0&&{1}&&&&\\
&&&{0}&&3&&0&&{1}&&&\\
&&{3}&&{0}&&6&&0&&{1}&&\\
&{0}&&15&&{0}&&10&&0&&{1}&\\
\bullet&&\bullet&&\bullet&&\bullet&&\bullet&&\bullet&&\bullet \end{array}
$$

\noindent which supplies in a compact  form the full defining information
on the triad.

This Hermite  triangle might be obtained  while considering subsequent powers of
$$
X =\left( \begin{array}{cccccc} 0&1&0&0&0&\ldots\\
1&0&1&0&0&\ldots\\
0&2&0&1&0&\ldots\\
0&0&3&0&1&\ldots\\
\ldots&\ldots&\ldots&\ldots&\ldots&\ldots  \end{array} \right)
$$

\noindent and then calculating  its  subsequent  rows $C_{n} = C_{0}
X^{n};\ n \geq 0$.

Let us take the opportunity to note that  Hermite polynomials  are known
to be the so called ``associated polynomials''  to    priority queue
organization of data bases  (see Propositions  4  and  6  in [28]).

\subsection{Laguerre  triad  and   Lah triangle}

It well known that
\begin{equation}\label{eq***5}
x^{n}=\sum\limits_{1\leq k\leq n} \left( n-1\\k-1 \right)\frac{n!}{k!}
(-1)^{k} L_{k}(x),\quad n,k\geq 0 ,
\end{equation}

\noindent where  $\{ L_{k}(x)\}_{k\geq 0}$  represents  self-inverse [11]
Laguerre binomial polynomial sequence which at the same time is  the unique
basic sequence of the  delta operator   $L=D/(D-1)$ [11, 12]  (see Chapter
4, Example 5 in [17] for extensions).  It is also well known that  the
following choice of numbers:
$$
q_{k}=2k;\ k\geq 0,\quad d_{k}=-k(k-1);\ k \geq 1,\quad i_{k}=-1;\ k \geq 0 ,
$$

\noindent  uniquely determines the following  recurrence for basic,
binomial  Laguerre polynomials
\begin{eqnarray}
&&xL_{k}(x)=-L_{k+1}(x)+2kL_{k}(x)-k(k-1)L_{k-1}(x),\quad k\geq 0 , \nonumber\\
&& \label{eq**5}\\
&&L_{0}(x)=1,\quad L_{-1}(x)=0.\nonumber
\end{eqnarray}

\noindent Connection constants  from  (\ref{eq***5})  are expressed by Lah
numbers [31] :
$$
L_{n,k}=(-1)^{k}C_{nk}=\left(\begin{array}{c} n-1\\n-k \end{array} \right)
\frac{n!}{k!}=\left( \begin{array}{c} n\\k \end{array}
\right)\frac{(n-1)!}{(k-1)!}. $$

\noindent Recall:  $L_{n,k}x^{k}$ number of  linear  deployments   of   $n$
distinguishable objects  in  exactly  $k$   cells  out of  $x\in  N$
distinguishable cells are available.  We infer from  (\ref{eq**5})  that
$c_{nk}$ satisfy accordingly
\begin{eqnarray}
&&c_{n+1,k}=-c_{n,k-1}+2kc_{n,k}-k(k+1)c_{n,k+1}, \nonumber\\
&&\label{eq*5}\\
&&c_{0,0}=1;\quad c_{0,k}=0 \quad {\rm for}\quad k>0.\nonumber
\end{eqnarray}

\noindent Hence for Ivo Lah numbers  we have:
\begin{eqnarray}
&&L_{n+1,k}=L_{n,k-1}+2kL_{n,k}+k(k+1)L_{n,k+1}, \nonumber \\
&&\label{eqL5}\\
&&L_{0,0}=1;\quad L_{0,k}=0 \quad {\rm for}\quad k>0.\nonumber
\end{eqnarray}

\noindent The transition matrix is then of the form
$$
X =\left( \begin{array}{cccccc} 0&-1&0&0&0&\ldots\\
0&2&-1&0&0&\ldots\\
0&2&4&-1&0&\ldots\\
0&0&-6&6&-1&\ldots\\
\ldots&\ldots&\ldots&\ldots&\ldots&\ldots  \end{array} \right),
$$

\noindent and the Laguerre triangle is obtained via  calculating its
subsequent rows  $C_{n} = C_{0} X^{n}$; $n \geq 0$. Thus we get:  Laguerre triangle
$$
\begin{array}{ccccccccccc}
&&&&&{1}&&&&&\\
&&&&{0}&&{-1}&&&&\\
&&&{0}&&-2&&{1}&&&\\
&&{0}&&-2&&6&&{-1}&&\\
&{0}&&{8}&&32&&-12&&{1}&\\
\bullet&&\bullet&&\bullet&&\bullet&&\bullet&&\bullet \end{array}
$$

\noindent and  Lah triangle
$$
\begin{array}{ccccccccccc}
&&&&&{1}&&&&&\\
&&&&{0}&&{1}&&&&\\
&&&{0}&&2&&{1}&&&\\
&&{0}&&4&&6&&{1}&&\\
&{0}&&{10}&&34&&12&&{1}&\\
\bullet&&\bullet&&\bullet&&\bullet&&\bullet&&\bullet \end{array}
$$

\noindent  with recurrence for  the corresponding triad polynomials being
of the form
\begin{eqnarray}
&&x\Lambda _{k}(x)=\Lambda_{k+1}(x)+2\Lambda_{k}(x)+k(k-1)\Lambda
_{k-1}(x),\quad k \geq 0, \nonumber \\
&&\label{eqLambda}\\
&&\Lambda _{0}(x)=1,\quad \Lambda _{-1}(x)=0.\nonumber
\end{eqnarray}

\subsection{Tchebychev  triad}

With  the following choice of    numbers  [1, 29]:
$$
q_{k}=0;\ k\geq 0,\quad d_{k}=1;\ k\geq 1,\quad i_{k}=1;\ k\geq 0 ,
$$

\noindent  the recurrence for corresponding  triad polynomials    takes the form
\begin{eqnarray}
&&x \Omega _{k}(x)=\Omega_{k-1}(x)+\Omega _{k+1}(x),\quad k \geq 0, \nonumber
\\
&&\label{eq**6} \\
&&\Omega _{0}=1,\quad \Omega _{-1} =0. \nonumber
\end{eqnarray}

\noindent which  is the recurrent definition of  $\Omega _{k}=U_{k}(x/2)$
polynomials where $U_{k}(x)$  represent classical  orthogonal Tchebychev
polynomials of the second kind. The Tchebychev triangle might be obtained
while considering subsequent powers of
$$
X =\left( \begin{array}{cccccc} 0&1&0&0&0&\ldots\\
1&0&1&0&0&\ldots\\
0&1&0&1&0&\ldots\\
0&0&1&0&1&\ldots\\
\ldots&\ldots&\ldots&\ldots&\ldots&\ldots  \end{array} \right),
$$

\noindent and then calculating  its   rows $C_{n} = C_{0} X^{n};\ n \geq
0$; thus getting the Tchebychev triangle
$$
\begin{array}{ccccccccccc}
&&&&&{1}&&&&&\\
&&&&{0}&&{1}&&&&\\
&&&{1}&&0&&{1}&&&\\
&&{0}&&2&&0&&{1}&&\\
&{2}&&{0}&&3&&0&&{1}&\\
\bullet&&\bullet&&\bullet&&\bullet&&\bullet&&\bullet \end{array}
$$

The connection constants   $c_{n,k}$ from
\begin{equation}\label{eq***6}
x^{n}=\sum\limits_{k\geq 0}c_{n,k}\Omega_{k}(x),\quad n \geq 0,
\end{equation}

\noindent are then determined  by   the  recurrence  dual to  (\ref{eq**6}):
\begin{eqnarray}
&&c_{n+1,k}=c_{n,k-1}+c_{n,k+1}, \nonumber\\
&&\label{eq*6}\\
&&c_{0,0}=1;\quad c_{0,k}=0\quad {\rm for}\quad k>0.\nonumber
\end{eqnarray}






\vspace{0.5cm}

\noindent {\erm Institute of Computer Science}

\noindent {\erm Bia\l ystok University }

\noindent{\erm  Sosnowa 64, PL-15-887 Bia\l ystok}

\noindent {\erm Poland}

\vspace{0.5cm}

\noindent Presented by Julian \L awrynowicz at the Session of the
Mathematical-Physical Commission of the \L \'od\'z Society of Sciences and
Arts on November 19, 2003

\vspace{0.5cm}
\noindent {\bf O TRIADACH DUALNYCH}
\vspace{0.2cm}

\noindent {\small S t r e s z c z e n i e}

{\small W pracy -- z inspiracji  opisu  dynamicznych baz danych w modelach
typu ``random walk'' [1, 2] -- wprowadza si\c{e} poj\c{e}cie  tzw. triad
dualnych. S\c{a} to uk\l ady dwu rekurencji dualnych dope\l nione
trzeci\c{a} relacj\c{a} o zazwyczaj wa\.znej interpretacji kombinatorycznej.
Podano szereg przyk\l ad\'ow takich triad ilustruj\c{a}c informacj\c{e} o
nich odpowiednimi ``tr\'oj\-k\c{a}tami''  na podobie\'nstwo  tr\'ojk\c{a}ta
Pascala. Triady traktowane s\c{a} jednocze\'snie jako uk\l ady dynamiczne z
czasem dyskretnym. }



\end{document}